\documentclass[11pt]{article}
\usepackage{amsmath,enumerate,amsfonts,color,amssymb,amsthm, verbatim}

\usepackage{hyperref}
\usepackage[normalem]{ulem}

\usepackage{tikz}

\def\RR{{\mathbb R}}
\def\Sphere{{\mathbb S}}
\def\mcA{{\mycal A}}

\newcommand{\f}{\varphi}

\def\bS{{\mathbb S}}
\def\eps{{\varepsilon}}

\newtheorem{theorem} {\sc  Theorem\rm}

\newtheorem{remark}[theorem]{\sc  Remark\rm}

\newcounter{marnote}

\DeclareFontFamily{OT1}{rsfs}{}
\DeclareFontShape{OT1}{rsfs}{m}{n}{ <-7> rsfs5 <7-10> rsfs7 <10-> rsfs10}{}
\DeclareMathAlphabet{\mycal}{OT1}{rsfs}{m}{n}

\def\be{\begin{equation}}
\def\ee{\end{equation}}

\newcommand{\R}{\mathbb{R}}

\def\be{\begin{equation}}
\def\ee{\end{equation}}
\def\bea#1\eea{\begin{align}#1\end{align}}

\newcommand{\blue}{\color{blue}}

\begin{document}
\title{Uniqueness of degree-one Ginzburg-Landau vortex in the unit ball in dimensions $N\geq 7$}

\author{Radu Ignat\thanks{Institut de Math\'ematiques de Toulouse \& Institut Universitaire de France, UMR 5219, Universit\'e de Toulouse, CNRS, UPS
IMT, F-31062 Toulouse Cedex 9, France. Email: Radu.Ignat@math.univ-toulouse.fr
}~, Luc Nguyen\thanks{Mathematical Institute and St Edmund Hall, University of Oxford, Andrew Wiles Building, Radcliffe Observatory Quarter, Woodstock Road, Oxford OX2 6GG, United Kingdom. Email: luc.nguyen@maths.ox.ac.uk}~, Valeriy Slastikov\thanks{School of Mathematics, University of Bristol, University Walk, Bristol, BS8 1TW, United Kingdom. Email: Valeriy.Slastikov@bristol.ac.uk}~ and Arghir Zarnescu\thanks{IKERBASQUE, Basque Foundation for Science, Maria Diaz de Haro 3,
48013, Bilbao, Bizkaia, Spain.}\,  \thanks{BCAM,  Basque  Center  for  Applied  Mathematics,  Mazarredo  14,  E48009  Bilbao,  Bizkaia,  Spain.
(azarnescu@bcamath.org)}\, \thanks{``Simion Stoilow" Institute of the Romanian Academy, 21 Calea Grivi\c{t}ei, 010702 Bucharest, Romania.}}

\date{}

\maketitle
\begin{abstract}
For $\eps>0$, we consider the Ginzburg-Landau functional for $\R^N$-valued maps defined in the unit ball $B^N\subset \R^N$ with the vortex boundary data $x$ on $\partial B^N$. 
In dimensions $N\geq 7$, we prove that for every $\eps>0$, there exists a unique global minimizer $u_\eps$ of this problem; moreover, $u_\eps$ is symmetric {and of the form}  
$u_\eps(x)=f_\eps(|x|)\frac{x}{|x|}$ for $x\in B^N$.  \\


\end{abstract}

\section{Introduction and main results}

In this note, we consider the following Ginzburg-Landau type energy functional
$$
E_\eps(u)=\int_{B^N} \Big[\frac{1}{2}|\nabla u|^2+\frac{1}{2\eps^2}W(1 - |u|^2)\Big]\,dx,
$$ 
where $\eps>0$, $B^N$ is the unit ball in $\R^N$, {$N \geq 2$}, and the potential $W\in C^1((-\infty,1] {\blue ;} \R)$ satisfies 
\begin{equation}
 W(0)=0,\, W(t)>0 \hbox{ for all } t\in (-\infty, 1]\setminus \{0\}, {\text{ and }} W \textrm{ is convex}.
 	\label{Eq:26VI18-E1}
\end{equation}
We investigate the global minimi{z}ers of the energy $E_\eps$ { in the set }
$$
\mcA:=\{ u\in H^1(B^N; \R^N):\,  u(x)=x \textrm{ on } \partial B^N = \bS^{N-1}\}.
$$ 
{The requirement that $u(x) = x$ on $\bS^{N-1}$ is sometimes referred {in the literature} as the vortex boundary condition.}

We note that in our analysis the convexity of $W$ needs not be strict; compare \cite{INSZ_AnnENS} where strict convexity is assumed.

The direct method in the calculus of variations yields the existence of a global minimizer $u_\eps$ of $E_\eps$ over $\mcA$ for all range of $\eps>0$. Moreover, any minimizer $u_\eps$ belongs to $C^1(\overline{B^N}; \R^N)$ and satisfies $|u_\eps| \leq 1$ and the system of PDEs (in the sense of distributions)
\be
\label{E-L}
-\Delta u_\eps=\frac1{\eps^2} u_\eps \, W'(1-|u_\eps|^2)\quad \textrm{ in } \, B^N.
\ee

The goal of this note is to give a short proof of the uniqueness and symmetry of the global minimizer of $E_\eps$ in $\mcA$ for \emph{all} $\eps>0$ in dimensions $N \geq 7$. We prove that,  in these dimensions, the global minimizer is unique and given by the unique radially symmetric critical point of $E_\eps$ defined by
\be
\label{def_sol_equa}
u_\eps(x)=f_\eps(|x|)\frac{x}{|x|} \quad \textrm{ for all } x\in B^N,
\ee
where the radial profile $f_\eps : [0,1]\to \R_+$ is the unique solution of 
\be
\label{1}
\left\{\begin{array}{l}
-f''_\eps-\frac{N-1}{r}f'_\eps+\frac{N-1}{r^2}f_\eps=\frac{1}{\eps^2} f_\eps \,W'(1-f_\eps^2) \quad 
\text{ for } r \in (0,1),\\
f_\eps(0)=0, f_\eps(1)=1.
\end{array}
\right.
\ee
Moreover, $f_\eps>0$ and $f'_\eps>0$ in $(0,1)$ (see e.g. \cite{ODE_INSZ}).

\medskip

\begin{theorem}
\label{GL_equator}
Assume that $W$ satisfies \eqref{Eq:26VI18-E1}. If $N\geq 7$, then for every $\eps>0$, $u_\eps$ given in \eqref{def_sol_equa} is the unique global minimizer of $E_\eps$  in $\mcA$.
\end{theorem}

To our knowledge, the question about the uniqueness of minimizers/critical points of $E_\eps$ in $\mcA$ for any $\eps > 0$ was raised in dimension $N = 2$ in the book of Bethuel, Brezis and H\'elein \cite[Problem 10, page 139]{BBH_book}, and in general dimensions $N \geq 2$ and also for the  blow-up limiting problem around the vortex (when the domain is the whole space $\R^N$ and by rescaling, $\eps$ can be assumed equal to $1$) in an article of Brezis \cite[Section 2]{Brezis99_PSPM}. 

It is well known that uniqueness is present for large enough $\eps>0$ for any $N \geq 2$. Indeed,  for any $\eps > (W'(1)/\lambda_1)^{1/2}$ where $\lambda_1$ is the first eigenvalue of $-\Delta$ in $B^N$ with  zero Dirichlet boundary condition, $E_\eps$ is strictly convex in $\mcA$ and thus has a unique critical point in $\mcA$ {(that is the global minimizer of our problem)}.

For \emph{sufficiently small} $\eps> 0$ all results regarding uniqueness question available in the literature are in the affirmative.  In particular, we have:
\begin{enumerate}[(i)]
\item Pacard and Rivi\`ere \cite[Theorem 10.2]{Pacard_Riviere} showed in dimension $N = 2$ that, for small $\eps > 0$, $E_\eps$ has in fact a unique critical point in $\mcA$.

\item Mironescu \cite{Mironescu_symmetry} showed in dimension $N = 2$ that, when $B^2$ is replaced by $\RR^2$ and $\eps = 1$, a local minimizer of $E_\eps$ subjected to a degree-one boundary condition at infinity is unique (up to translation and suitable rotation). This was generalized to dimension $N = 3$ by Millot and Pisante \cite{mil-pis} and dimensions $N \geq 4$ by Pisante \cite{Pisante11-JFA}, also in the case of the blow-up limiting problem on $\RR^N$ and $\eps = 1$.
\end{enumerate}
These results should be compared to those for the limit problem on the unit ball obtained by sending $\eps \rightarrow 0$. In this limit, the Ginzburg-Landau problem `converges' to the harmonic map problem from $B^N$ to $\bS^{N-1}$. It is well known that, the vortex boundary condition gives rise to a unique minimizing harmonic map $x \mapsto \frac{x}{|x|}$ if $N\geq 3$; see Brezis, Coron and 
Lieb~\cite{BrezisCoronLieb}  in dimension $N = 3$, J\"{a}ger and Kaul \cite{JagKaul} in dimensions $N \geq 7$, and Lin \cite{Lin-CR87} in dimensions $N \geq 3$.

We highlight that, in contrast to the above, our result holds for \emph{all} $\eps > 0$, provided that $N \geq 7$. The method of our proof deviates somewhat from that in the aforementioned works. In fact it is reminiscent of our recent work \cite{INSZ_AnnENS} on the (non-)uniqueness and symmetry of minimizers of the Ginzburg-Landau functionals for $\RR^M$-valued maps defined on $N$-dimensional domains, where $M$ is not necessarily the same as $N$. However we note that the results in \cite{INSZ_AnnENS} do not  directly apply to the present context, as in \cite{INSZ_AnnENS} it is required that $W$ be \emph{strictly convex}. Furthermore, a priori, it is not clear why non-strict convexity of the potential $W$ is sufficient to ensure uniqueness of global minimizers.

We exploit the convexity of $W$ to lower estimate the `excess' energy by a suitable quadratic energy which can be handled by the factorization trick \`a la Hardy. Indeed, the positivity of the excess energy is then related to the validity of a Hardy-type inequality, which explains our restriction of $N \geq 7$. This echoes our observation made in \cite{INSZ_AnnENS} that a result of J\"ager and Kaul \cite{JagKaul} on the minimality of the equator map in these dimensions is related to a certain inequality involving the sharp constant in the Hardy inequality.

We expect that our result remains valid in dimensions $2 \leq N \leq 6$, but this goes beyond the scope of this note and remains for further investigation.

\section{Proof of Theorem \ref{GL_equator}}

Theorem \ref{GL_equator} will be obtained as a consequence of a stronger result on the uniqueness of global minimizers of for the $\RR^M$-valued Ginzburg-Landau functional with $M \geq N$. By a slight abuse of notation, we consider the energy functional
$$
E_\eps(u)=\int_{B^N} \Big[\frac{1}{2}|\nabla u|^2+\frac{1}{2\eps^2}W(1 - |u|^2)\Big]\,dx,
$$ 
where $u$ belongs to 
$$
\mcA:=\{ u\in H^1(B^N; \R^M):\,  u(x)=x \textrm{ on } \partial B^N = \bS^{N-1} \subset \RR^M\}.
$$ 

\begin{theorem}
\label{GL_equatorNM}
Assume that $W$ satisfies \eqref{Eq:26VI18-E1}. If $M \geq N\geq 7$, then for every $\eps>0$, $u_\eps$ given in \eqref{def_sol_equa} is the unique global minimizer of $E_\eps$  in $\mcA$.
\end{theorem}

When $W$ is strictly convex, the above theorem is proved in \cite{INSZ_AnnENS}; see Theorem 1.7. The argument therein uses the strict convexity in a crucial way.

\begin{proof}

The proof will be done in several steps. First, we consider the difference between the energies of the critical point $u_\eps$,  defined in \eqref{def_sol_equa},  and an arbitrary competitor $u_\eps+v$ and show that this difference is controlled from below by some quadratic energy functional $F_\eps(v)$. Second, we employ the positivity of the radial profile $f_\eps$ in \eqref{1} and apply the Hardy decomposition method in order to show that $F_\eps(v)\geq 0$, which proves in particular that $u_\eps$ is a global minimizer of $E_\eps$. Finally, we characterise the situation when this difference is zero and conclude to the uniqueness of the global minimizer $u_\eps$. 

\bigskip
\par\noindent{\it Step 1: Lower bound for energy difference.} 
For any  $v\in H^1_0(B^N; \R^M)$, we have
\begin{align*}
E_\eps(u_\eps+v)-E_\eps(u_\eps)
	&=\int_{B^N}\Big[ \nabla u_\eps\cdot\nabla v +\frac{1}{2}|\nabla v|^2\Big]\,dx\\
		&\qquad +\frac{1}{2\eps^2}\int_{B^N} \Big[W(1-|u_\eps+v|^2)-W(1-|u_\eps|^2)\Big]\,dx.
\end{align*}
Using the convexity of $W$, we have
\[
W(1-|u_\eps+v|^2)-W(1-|u_\eps|^2)\ge -W'(1-|u_\eps|^2)(|u_\eps+v|^2-|u_\eps|^2).
\]
The last two relations imply that
\begin{align*}
E_\eps(u_\eps+v)-E_\eps(u_\eps)
	&\ge \int_{B^N} \Big[\nabla u_\eps\cdot\nabla v-\frac{1}{\eps^2} W'(1-f_\eps^2)u_\eps\cdot v\Big]\,dx\nonumber\\
		&\qquad\qquad +\int_{B^N} \Big[ \frac{1}{2} |\nabla v|^2-\frac{1}{2\eps^2}W'(1-f_\eps^2)|v|^2\Big]dx.
\end{align*}
Moreover, by \eqref{E-L}, we obtain
\bea\label{rel:endif3}
E_\eps(u_\eps+v)-E_\eps(u_\eps) \ge  \int_{B^N} \Big[ \frac 12 |\nabla v|^2-\frac{1}{2\eps^2} W'(1-f_\eps^2) |v|^2\Big]\,dx =: \frac12 F_\eps(v)
\eea
for all  $v\in H^1_0({B^N};\R^M)$.

\bigskip
\par\noindent{\it Step 2: A rewriting of $F_\eps(v)$ using the decomposition $v=f_\eps w$ for every scalar test function $v\in C^\infty_c(B^N\setminus \{0\}; \R)$.} 
We consider the operator
$$L_\eps:=\frac12 \nabla_{L^2}F_\eps=-\Delta-\frac1{\eps^2} W'(1-f_\eps^2).$$
Using the decomposition $$v=f_\eps w$$ for the scalar function $v\in C^\infty_c(B^N\setminus \{0\}; \R)$, we have (see e.g. \cite[Lemma~A.1]{INSZ3}):
\begin{align*}
F_\eps(v)=\int_{B^N} L_\eps v\cdot v\, dx&=\int_{B^N} w^2 L_\eps f_\eps\cdot f_\eps\, dx+\int_{B^N} f_\eps^2 |\nabla w|^2\, dx\\
& =\int_{B^N} f_\eps^2 \bigg( |\nabla w|^2-\frac{N-1}{r^2}w^2\bigg)\, dx,
\end{align*}
because \eqref{1} yields $L_\eps f_\eps\cdot f_\eps=-\frac{N-1}{r^2} f_\eps^2$ in ${B^N}$. 

\bigskip
\par\noindent{\it Step 3: We prove that  $F_\eps(v)\geq 0$ for every scalar test function $v\in C^\infty_c(B^N\setminus \{0\}; \R)$.} 
Within the notation $v=f_\eps w$ of Step 2 with $v,w \in  C^\infty_c(B^N\setminus \{0\}; \R)$, we use the decomposition $$w=\f g$$ with $\f=|x|^{-\frac{N-2}{2}}$ being the first eigenfunction of the Hardy's operator $-\Delta-\frac{(N-2)^2}{4|x|^2}$ in $\R^N\setminus\{0\}$ and $g\in C^\infty_c({B^N}\setminus \{0\}; \R)$. We compute
\begin{align*}
|\nabla w|^2=|\nabla \f|^2 g^2+|\nabla g|^2 \f^2+\frac12 \nabla (\f^2)\cdot \nabla (g^2).
\end{align*}
As $|\nabla \f|^2=\frac{(N-2)^2}{4|x|^2}\f^2$ and $\f^2$ is harmonic in ${B^N}\setminus\{0\}$, integration by parts yields
\begin{align}
\nonumber
F_\eps(v)&=\int_{B^N} f_\eps^2 \bigg( |\nabla g|^2 \f^2+\frac{(N-2)^2}{4r^2}\f^2g^2-\frac{N-1}{r^2}\f^2 g^2\bigg)\, dx-\frac12 \int_{{B^N}} \nabla (\f^2)\cdot \nabla (f_\eps^2)g^2\, dx\\
\nonumber &\geq \int_{B^N} f_\eps^2 |\nabla g|^2 \f^2\, dx+\bigg(\frac{(N-2)^2}{4}-(N-1)\bigg)\int_{{B^N}}\frac{f_\eps^2}{r^2}\f^2g^2\, dx\\
\label{23}
&\geq \bigg(\frac{(N-2)^2}{4}-(N-1)\bigg)\int_{{B^N}}\frac{v^2}{r^2}\, dx\geq 0,
\end{align}
where we have used $N\geq 7$ and $\frac12\nabla (\f^2)\cdot \nabla (f_\eps^2)=2\f \f' f_\eps f'_\eps\leq 0$ in ${B^N}\setminus \{0\}$.

\bigskip
\par\noindent{\it Step 4: We prove that $F_\eps(v)\geq 0$ for every $v\in H^1_0(B^N; \R^M)$ meaning that $u_\eps$ is a global minimizer of 
$E_\eps$ over $\mcA$; moreover, $F_\eps(v)=0$ if and only if $v=0$.} Let $v\in H^1_0(B^N; \R^M)$. As a point has zero $H^1$ capacity in $\R^N$, a standard density argument implies the existence of a sequence $v_k\in C^\infty_c({B^N}\setminus \{0\}; \R^M)$ such that $v_k\to v$ in $H^1(B^N, \R^M)$ and a.e. in $B^N$. On the one hand, by 
definition \eqref{rel:endif3} of $F_\eps$, since $W'(1-f_\eps^2)\in L^\infty$, we deduce that $F_\eps(v_k)\to F_\eps(v)$ as $k\to \infty$.
On the other hand, by \eqref{23} and Fatou's lemma, we deduce
\begin{align*}
\liminf_{k\to \infty} F_\eps(v_k)&\geq  \bigg(\frac{(N-2)^2}{4}-(N-1)\bigg)\liminf_{k\to \infty} \int_{{B^N}}\frac{v_k^2}{r^2}\, dx\\
&\geq \bigg(\frac{(N-2)^2}{4}-(N-1)\bigg) \int_{{B^N}}\frac{v^2}{r^2}\, dx.
\end{align*}
Therefore, we conclude that
$$F_\eps(v)\geq \bigg(\frac{(N-2)^2}{4}-(N-1)\bigg) \int_{{B^N}}\frac{v^2}{r^2}\, dx\geq 0,\quad \forall v\in H^1_0(B^N; \R^M),$$
implying by \eqref{rel:endif3} that $u_\eps$ is a minimizer of $E_\eps$ over $\mcA$.
Moreover, $F_\eps(v)=0$ if and only if $v=0$.

\bigskip
\par\noindent{\it Step 5: Conclusion.} We have shown that $u_\eps$ is a global minimizer. Assume that $\tilde u_\eps$ is  another global minimizer of $E_\eps$ over $\mcA$. If $v:=\tilde u_\eps-u_\eps$,
then $v\in H^1_0(B^N; \R^M)$ and by Steps 1 and 4, we have that $0=E_\eps(\tilde u_\eps)-E_\eps(u_\eps)\geq F_\eps(v)\geq 0$, which yields $F_\eps(v)=0$. Step 4 implies that $v=0$, i.e., $\tilde u_\eps=u_\eps$.
 \end{proof}

\begin{remark}
Recall that in the case $M \geq N \geq 7$, J\"{a}ger and Kaul \cite{JagKaul} proved the uniqueness of global minimizer for harmonic map problem
$$
\min_{u \in \mcA_*} \int_{B^N} |\nabla u|^2 \, dx,
$$
where $\mcA_* = \{  u\in H^1(B^N; \Sphere^{M-1}):\,  u(x)=x \textrm{ on } {\partial B^N} = \bS^{N-1} \subset \bS^{M-1} \}$. This can also be seen by the method above as observed in our earlier paper \cite{INSZ_AnnENS}. We give the argument here for readers' convenience: Take a perturbation $v \in H_0^1(B^N, \R^M)$ of the harmonic map $u_*(x) = \frac{x}{|x|}$ such that $|u_*(x) + v(x)|=1$ a.e. in $B^N$. Then, by \cite[Proof of Theorem~5.1]{INSZ_AnnENS},
\[
\int_{B^N} \Big[|\nabla (u_*+v)|^2 - |\nabla u_*|^2\Big]\,dx =\int_{B^N} \Big[|\nabla v|^2- |\nabla u_*|^2|v|^2\Big]\,dx= \int_{B^N}\Big[|\nabla v|^2- (N-1)\frac{|v|^2}{|x|^2}\Big]\,dx.
\]
Using Hardy's inequality in dimension $N$ we arrive at 
$$
\int_{B^N} \Big[|\nabla (u_*+v)|^2 - |\nabla u_*|^2\Big]\,dx \geq \left(\frac{(N-2)^2}{4} - (N-1)\right) \int_{B^N} \frac{|v|^2}{|x|^2}\,dx.
$$
The result follows since $N \geq 7$.
\end{remark}

\section*{Acknowledgment.} R.I. acknowledges partial support by the ANR project ANR-14-CE25-0009-01. V.S. acknowledges  support by  the Leverhulme
grant RPG-2014-226. A.Z. was partially supported by a Grant of the Romanian National
Authority for Scientific Research and Innovation, CNCS-UEFISCDI, project number
PN-II-RU-TE-2014-4-0657; by the Basque Government through the BERC
2014-2017 program; and by the Spanish Ministry of Economy and Competitiveness
MINECO: BCAM Severo Ochoa accreditation SEV-2013-0323.

\bibliographystyle{acm}
\bibliography{paris.bib,LiquidCrystals.bib}

\def\cprime{$'$}
\begin{thebibliography}{10}

\bibitem{BBH_book}
{\sc Bethuel, F., Brezis, H., and H{\'e}lein, F.}
\newblock {\em Ginzburg-{L}andau vortices}.
\newblock Progress in Nonlinear Differential Equations and their Applications,
  13. Birkh\"auser Boston Inc., Boston, MA, 1994.

\bibitem{Brezis99_PSPM}
{\sc Brezis, H.}
\newblock Symmetry in nonlinear {PDE}'s.
\newblock In {\em Differential equations: {L}a {P}ietra 1996 ({F}lorence)},
  vol.~65 of {\em Proc. Sympos. Pure Math.} Amer. Math. Soc., Providence, RI,
  1999, pp.~1--12.

\bibitem{BrezisCoronLieb}
{\sc Brezis, H., Coron, J.-M., and Lieb, E.~H.}
\newblock Harmonic maps with defects.
\newblock {\em Comm. Math. Phys. 107}, 4 (1986), 649--705.

\bibitem{ODE_INSZ}
{\sc Ignat, R., Nguyen, L., Slastikov, V., and Zarnescu, A.}
\newblock Uniqueness results for an {ODE} related to a generalized
  {G}inzburg-{L}andau model for liquid crystals.
\newblock {\em SIAM J. Math. Anal. 46}, 5 (2014), 3390--3425.

\bibitem{INSZ3}
{\sc Ignat, R., Nguyen, L., Slastikov, V., and Zarnescu, A.}
\newblock Stability of the melting hedgehog in the {L}andau-de {G}ennes theory
  of nematic liquid crystals.
\newblock {\em Arch. Ration. Mech. Anal. 215}, 2 (2015), 633--673.

\bibitem{INSZ_AnnENS}
{\sc Ignat, R., Nguyen, L., Slastikov, V., and Zarnescu, A.}
\newblock On the uniqueness of minimisers of {G}inzburg-{L}andau functionals.
\newblock {\em arXiv:1708.05040\/} (2017).

\bibitem{JagKaul}
{\sc J\"ager, W., and Kaul, H.}
\newblock Rotationally symmetric harmonic maps from a ball into a sphere and
  the regularity problem for weak solutions of elliptic systems.
\newblock {\em J. Reine Angew. Math. 343\/} (1983), 146--161.

\bibitem{Lin-CR87}
{\sc Lin, F.-H.}
\newblock A remark on the map {$x/|x|$}.
\newblock {\em C. R. Acad. Sci. Paris S\'er. I Math. 305}, 12 (1987), 529--531.

\bibitem{mil-pis}
{\sc Millot, V., and Pisante, A.}
\newblock Symmetry of local minimizers for the three-dimensional
  {G}inzburg-{L}andau functional.
\newblock {\em J. Eur. Math. Soc. (JEMS) 12}, 5 (2010), 1069--1096.

\bibitem{Mironescu_symmetry}
{\sc Mironescu, P.}
\newblock Les minimiseurs locaux pour l'\'equation de {G}inzburg-{L}andau sont
  \`a sym\'etrie radiale.
\newblock {\em C. R. Acad. Sci. Paris S\'er. I Math. 323}, 6 (1996), 593--598.

\bibitem{Pacard_Riviere}
{\sc Pacard, F., and Rivi\`ere, T.}
\newblock {\em Linear and nonlinear aspects of vortices}, vol.~39 of {\em
  Progress in Nonlinear Differential Equations and their Applications}.
\newblock Birkh\"auser Boston, Inc., Boston, MA, 2000.
\newblock The Ginzburg-Landau model.

\bibitem{Pisante11-JFA}
{\sc Pisante, A.}
\newblock Two results on the equivariant {G}inzburg-{L}andau vortex in
  arbitrary dimension.
\newblock {\em J. Funct. Anal. 260}, 3 (2011), 892--905.

\end{thebibliography}

\end{document}